\documentclass[10pt]{amsart}
\usepackage{amssymb, amscd, amsmath, amsthm, epsf, epsfig, latexsym}
\def\zz{{\bf Z}}
\def\ff{{\bf F}}
\def\qq{{\bf Q}}

\def\rr{{\bf R}}
\def\dd{{\bf D}}
\def\co{\colon\thinspace}
\def\calc{\mathcal{C}}

\def\calg{\mathcal{G}}

\def\calg{\mathcal{G}}

\def\co{\colon}
\newcommand\disc{\operatorname{disc}}
\newcommand\Disc{\operatorname{Disc}}
\newcommand\Res{\operatorname{Res}}

\newtheorem{theorem}{Theorem}[section]
\newtheorem{lemma}[theorem]{Lemma}
\newtheorem{corollary}[theorem]{Corollary}

\theoremstyle{definition}
\newtheorem{definition}[theorem]{Definition}
\numberwithin{equation}{section}

\begin{document}
\title{The Algebraic Concordance Order of a Knot}
\author{Charles Livingston}
\thanks{This work was supported by a grant  from the NSF}
\address{Department of Mathematics, Indiana University, Bloomington, IN 47405}
\email{livingst@indiana.edu}


\maketitle

Since the inception of knot concordance, questions related to torsion in the concordance group, $\calc_1$, have been of particular interest; see for instance~\cite{fox1, go2, k}. The only known torsion in $\calc_1$ is 2--torsion, arising from amphicheiral knots, whereas  Levine's analysis of higher dimensional concordance revealed far more 2--torsion and  also 4--torsion in $\calc_{2n-1}, n>1$.   Casson and Gordon~\cite{cg1, cg2} demonstrated that Levine's algebraic classification of concordance does not apply to $\calc_1$;  since then, the basic questions relating to torsion in $\calc_1$ have remained open.  However, many of the deep theoretical tools of 4--dimensional topology, for instance~\cite{cg1, do, os}, have been applied to this problem,  ruling out potential classes of order two and four.  Examples of this work includes~\cite{ grs, jn, lis, ln1, ln2, ta}.  

Surprisingly, the determination of the algebraic order of a knot, as defined by Levine, has remained a difficult problem.  The only work analyzing a set of examples was a fairly technical paper by Morita~\cite{mo}  studying knots with 10 or fewer crossings.  Our main goal   is to study algebraic concordance to the extent necessary to easily determine the concordance order of knots; we demonstrate the effectiveness of these results by determining the algebraic orders of all 2,977 prime knots of 12 or fewer crossings.

\vskip.1in
\noindent{\bf Background and Summary of Results}

In 1969 Levine~\cite{le1} defined the algebraic knot concordance groups, $\calg^\zz_{\pm 1}$,   defined   homomorphisms $\phi_n \co \calc_n \to \calg^\zz_{(-1)^n}$ with domain the concordance group
of knotted $(2n-1)$--spheres in $S^{2n+1}$,  and   proved that for $n >1$, $\phi_n$ is essentially an isomorphism. In a second paper,~\cite{le2}, he  gave a complete set of invariants that determine a class in $\calg^\zz_{\pm 1}$.  Consequently, he proved that $ \calg^\zz_{\pm 1} $ is isomorphic to the infinite direct sum $\oplus_\infty \zz \oplus_\infty \zz/2\zz \oplus_\infty\zz/4\zz$.     

Our interest is in the case $n=1$, so we now drop the subscripts from $\calg, \calc$ and $\phi$.  As described by Levine, $\calg^\zz $ injects into a {\it rational algebraic concordance group} $\calg_\qq$.  There are similarly defined groups over other fields, $\calg_\ff$.
A key result from~\cite{le2} is the following.

\vskip.1in
\noindent{\bf Theorem.} {\it A class in $\calg_\qq$ is trivial if and only if it is trivial in $\calg_\ff$ for $\ff = \rr$ and $\ff = \qq_p$ for all primes $p$,  where   $\qq_p$ denotes the $p$--adic rationals.}
\vskip.1in

We have four principal goals in this paper.    We first have a theoretical result  which implies that 
the classification of algebraic knot concordance is effectively computable by restricting the set of primes that need to be considered.  Recall that classes in $\calg^\zz$ are represented by Seifert matrices.

\vskip.1in
\noindent{\bf Theorem.} {\it For a nonsingular integral Seifert matrix $V$, the  class   $[V] \in \calg^\zz$ is of infinite order if and only if it is nontrivial in $\calg_\rr$; if it is of finite   order,  it is of order $4$ if and only if it is of order $4$ in $\calg_{\qq_p}$ for some $p$ dividing $\Delta_V(-1)$ with $p \equiv 3 \mod 4$;  if  it is of order $2$, then  it is of order  $2$ in $\calg_{\qq_p}$ for some prime $p$ dividing  $2\det(V) \Disc(\bar{\Delta}_V(t))$ }
\vskip.1in

\noindent Here $\Delta_V(t) =  \det(V - tV^t$)   is the Alexander polynomial,    $\bar{\Delta}_V(t)$ denotes the product of the irreducible factors of   $\Delta_V(t)$, and $\Disc$ denotes the discriminant of a polynomial, reviewed in Appendix~\ref{discsection}.  (The result for 4--torsion was essentially proved by Morita~\cite{mo}, though he did not eliminate the prime 2 or primes $p \equiv 1 \mod 4$.) 

This result, calling on an analysis of Witt groups over the $p$--adics, can be difficult to apply.  Thus, a second goal of this paper is to provide simplifying criteria. One example is the following.

\vskip.1in
\noindent{\bf Theorem.} {\it If $V$ is an integral Seifert matrix representing a class of order 4 in $\calg^\zz$, then for some $p  \equiv 3 \mod 4$ and some symmetric irreducible factor $g$ of $\Delta_V(t)$, $p$ divides $g(-1)$ and $g$ has odd exponent in $\Delta_V(t)$}.
\vskip.1in

A third goal of the paper is to apply the analysis to specific knots.  We do this by determining the algebraic order of all 2,977 prime knots of 12 or fewer crossings.

The final goal of the paper is to provide an expository account of some of the underlying theory that is not included in Levine and which might be unfamiliar to many geometric topologist.  Much of this background is provided in appendices concerning $p$--adic numbers, Witt groups, and discriminants.

   Notice that the results quoted above concern classes in $\calg^\zz$ and do not apply to the full rational group $\calg_\qq$.  A complete classification of $\calg^\zz$ was presented by Stoltzfus in~\cite{st}.  The results we present here, using Witt groups of quadratic forms over finite fields, could be also obtained via the number theoretic approach of~\cite{st}.     This paper can be viewed as an effort to most simply extract from the structure of the integral algebraic concordance group enough information to be able to effectively analyze torsion. 
\vskip.1in
\noindent{\bf Outline}  In the first section we review   Levine's algebraic concordance group.   Section 2 quickly gives the identification of  elements of infinite order, as was done in~\cite{le2}, and also notes that $\calg_\rr$ is torsion free.  In Section 4 we consider 4-torsion, extending Morita's result,~\cite{mo}, and developing simpler tests for detecting elements of order 4.  Torsion of order 2 is analyzed in Section 4.  In Section 5 we describe how the earlier work of this article can be applied to specific examples by analyzing all prime knots of 12 or fewer crossings.  Only one case calls on explicit work over the $p$--adics.   Appendices A, B, and C, provide the necessary background material on $p$--adic algebra, Witt groups, and the theory of polynomial discriminants and resultants.

\vskip.1in
\noindent{\it Acknowledgments} In learning the algebra  needed here, the assistance of  Jim Davis, Darrell Haile, and Michael Larsen was invaluable.  I also thank Neal Stoltzfus for helpful discussions.

\section{Review of Levine's Algebraic Concordance Group}

To each knot $K$ one can associate an integer Seifert matrix $V$; that is, an integer matrix of size $2g \times 2g$ for some $g$, satisfying $\det(V - V^t) = \pm1$.  A Seifert matrix is called {\it algebraically slice}, or {\it Witt trivial}, if there is a rank $g$ submodule of $\zz^{2g}$ on which the  bilinear form associated to $V$ vanishes. The abelian group $\calg^\zz$ is  the Witt group of  Seifert matrices, defined to be the set of equivalence classes of such $V$, where $V \sim W$ if $V \oplus -W$ is Witt trivial.  Addition is via direct sum.

Considering rational matrices instead, one can form the group $\calg^\qq$.  Here one restricts to $V$ satisfying $\det( (V - V^t)(V+V^t)) \ne 0$.  The inclusion $\calg^\zz \to \calg^\qq$ is injective~\cite[Section 3]{le2}.

Levine proved that every element in $\calg^\qq$ has an invertible representative, and in fact the same is true in $\calg^{\zz}$.  In brief, if $\det(V)= 0$, one can perform row operations over $\zz$ to form a matrix with bottom row identically 0.  Performing the corresponding column operations preserves this feature.  Further simultaneous row and column operations results in a matrix in the form

 $$\left(
\begin{array}{ccc}
  B&   b_1 & 0   \\
  b_2 & b_3   & 1  \\   
  0&  0   & 0  \\   
\end{array}
\right),
$$ and one then proves that $V$ is Witt equivalent to $B$. Given this, we assume henceforth that Seifert matrices are invertible.  In particular, the Alexander polynomial $\Delta_V(t)  = \det(V - tV^t)$ is of degree $2g$ with leading and constant coefficient $\det(V)$.

  Associated to each element $V \in \calg^\qq$ is a triple   $(M, Q, T)$, where $M$ is a $2g$--dimensional rational vector space, $Q = V + V^t$ is a nonsingular symmetric bilinear form on $M$, and $T$ is the linear transformation $V^{-1}V^t$, an  isometry of $(M,Q)$.  (Levine considered $T = -V^{-1}V^t$; we change signs to be consistent with the standard definition of the Alexander polynomial, $\Delta_V(t) = \det( V - tV^t)$; we denote the characteristic polynomial of $T$ by $\Delta_T(t) = \det (T - tI)$.)   We have  the following.  
  
\begin{theorem} If $V$ is a nonsingular Seifert matrix and $(M, Q, T)$ is the associated isometric structure, then $\Delta_T(t) = \det(V) \Delta_V(t)$ and  $\Delta_T(1) \Delta_T(-1)  \ne 0$.   
\end{theorem}
 
 \noindent{\bf Convention} {\it When working with elements $(M,Q,T)$ in an algebraic concordance group, we will necessarily work with $\Delta_T(t)$.  When presenting results that can be applied directly to knots, we will, when possible, work with $\Delta_V(t)$.}
 \vskip.1in

  The algebraic concordance group, $\calg_\qq$, is defined to be the group of Witt classes of such {\it isometric structures} $(M,Q,T)$, ($T$ is required to satisfy $\Delta_T(1)\Delta_T(-1) \ne 0$) where such a structure is Witt trivial if $M$ contains a $g$--dimensional subspace that is invariant under $T$ and on which $Q$ vanishes.  We have the following result~\cite[Theorem 8]{le2}.

\vskip.1in

\begin{theorem}  The map $\calg^\qq \to \calg_\qq$ is an isomorphism.
\end{theorem}

Given any field $\ff$ there is a similarly defined algebraic concordance group of isometric structures:  $\calg_\ff$.  We have the following~\cite[Proposition 17]{le2}.
\vskip.1in
\begin{theorem}  A class $(M, Q, T) \in \calg_\qq$ is trivial if and only if $(M,Q,T)$ represents $0 \in \calg_\ff$ for $\ff = \rr$ and $\ff = \qq_p$ for all $p$.\end{theorem}

\noindent  Here $\qq_p$ denotes the $p$--adic rationals.  See Section~\ref{padicsection} for a review of $p$--adic numbers.

For every  polynomial $f \in \ff[t]$ there is a Witt group of isometric structures $\calg_\ff^f$ defined as above, but restricting to those structures for which the characteristic polynomial of $T$ is a power of $f$.  According to~\cite[Lemma 11]{le2}:

\begin{theorem}  $\calg_\ff \cong \oplus_{\delta } \calg_\ff^\delta$ where the sum is over all irreducible symmetric polynomials.  In particular, a class in $(M,Q,T) \in \calg_\ff$ is trivial if and only if its projection $(M^\delta, Q^\delta, T^\delta)$ in  $\calg_\ff^\delta$ is trivial for all symmetric irreducible factors $\delta$ of the characteristic polynomial of $T$.  (Symmetric means that $\delta(t) = at^{k}\delta(t^{-1})$ for some integer $k$ and some field element $a$.)
\end{theorem}

\noindent We need to expand on this briefly.  Suppose that $\Delta_T(t) $ factors as $\prod_i \delta_i(t)^{k_i} \prod_j g_j(t)^{l_j}$, where the $\delta_i$ are distinct irreducible symmetric factors and the $g_i$ are the remaining irreducible factors. Since $\Delta_T(t)$ is monic, we can choose each factor to be monic.   Let $\hat{\delta}_i = \Delta_T / \delta_i^{k_i}$.  Then $M^{\delta_i} = \text{Im} (\hat{\delta}_i ^N(T))$ for any large $N$.  (This follows from the fact that $\hat{\delta}_i(T)^N$ is an isomorphism of $M^{\delta_i}$ and annihilates all the other $M^f$ summands if $N$ is large.)  In addition,  $T$ restricted to $M^{\delta_i}$ has characteristic polynomial $\delta_i^k(t)$ for some $k$, since $\delta_i$ was chosen to be monic.

We next have~\cite[Proposition 16]{le2}, stating the following.
 
\begin{theorem} \label{levine3}   A class $(M,Q,T) \in \calg_\ff^\delta $, where $\ff = \rr$ or $\ff = \qq_p$ and $\delta$ is irreducible symmetric, is trivial if and only if the characteristic polynomial of $T$, $\Delta_T(t)$, is $\delta^e$ with $e$ even and the form $(M,Q)$ is trivial in the Witt group of $\ff$, $W(\ff)$.
\end{theorem}

\noindent The Witt groups of symmetric bilinear forms $W(\ff)$ are reviewed in Appendix~\ref{wittgroupsection}.
According to \cite[Lemma 7]{le2},

\begin{theorem}   Let $(M, Q, T)$ be an isometric structure over a field $\ff$ (with the property that the characteristic polynomial of $T$, $\Delta_T(t)$, satisfies $\Delta_T(1) \Delta_T(-1) \ne 0$).  Then:  
\begin{enumerate}
\item $\Delta_T(t) = t^d \Delta_T(t^{-1})$ with
  $d = \deg(\Delta_T(t))$ even.
\item If $(M,Q,T) = 0 \in \calg_\ff$ then
  $\Delta_T(t)= ct^e  f(t)f(t^{-1})$, where      $f(t)$ is a polynomial of
   degree $e$, and $c \in \ff$.

\item $\det(Q) = \Delta_T(1) \Delta_T(-1) \in \ff^* / (\ff^*)^2$.
\end{enumerate} 
\end{theorem}

Note that since $d$ is even, $\Delta_T(t) = \det(T - tI)$ is monic with leading coefficient 1.  Thus, the direct sum decomposition $\calg_\ff \cong \oplus_{\delta } \calg_\ff^\delta$ can be taken over polynomials that are irreducible, symmetric, and monic.

\section{Elements of Infinite order}

The Witt groups $W(\qq_p)$ are all finite (see Appendix~\ref{wittgroupsection}).  It follows that any element of infinite order in $\calg_\qq$ is of infinite order in  $\calg_\rr$.  The class $(M,Q,T) \in \calg_\rr$ splits as the direct sum of classes in $\calg_{\rr}^\delta$ where $\delta$ is an irreducible symmetric real polynomial.  The only such polynomials are, up to a unit,  of the form   $t^2 +2at + 1$, where $a^2   < 1$.  The complex roots of this polynomial are the unit complex numbers $\omega$, where $\omega = e^{i \theta}$ and $\cos \theta = a$.

On $\calg_\rr^\delta$ we have the surjective signature function $\sigma \co \calg_\rr^\delta \to 2\zz \subset \zz$ defined by  $\sigma(M,Q, T)= \text{sign}(Q)$.

\begin{theorem}  The signature function $\sigma$ is an isomorphism. 
\end{theorem}
 \begin{proof} As described in Appendix~\ref{wittgroupsection},  the signature function defines an isomorphism of $W(\rr)$ with $2\zz$.  Thus, by Theorem~\ref{levine3}, a nontrivial form $(M,Q,T) \in \calg_\rr^\delta $ in the kernel of  the signature function on $ \calg_\rr^\delta$ would have signature is 0 and $\Delta_T(t) = \delta(t)^k$ for some odd $k$. This can be seen to be impossible as follows.  If $k$ is odd, then $M$ is of dimension $4m +2$.  The determinant of $Q$ is given, modulo squares, by $\delta(1)\delta(-1) = (2+2a)(2-2a) = 4(1 - a^2) >0$.  However, a diagonal form of dimension $4m + 2$ of signature 0 has determinant $-1$.
\end{proof} 
  
The following theorem provides a means of computing the associated signatures.

\begin{theorem} The signature function \text{\rm sign}$((1-\omega)V + (1 - \bar{\omega})V^t)$ defined for $\omega \in S^1$ has jumps only at the unit roots of the Alexander polynomial.  If   $\omega = e^{i \theta}$, with $\cos \theta = a$, is a root of   $\Delta_V(t)$, then $\delta_a(t) = t^2 +2at + 1$ is a factor of $\Delta_V(t)$ and the jump in the signature    \text{\rm sign}$((1-\omega)V + (1 - \bar{\omega})V^t)$ at $\omega$ equals, up to sign, the signature of $V+ V^t$ restricted to $\calg_\rr^{\delta_a}$.
\end{theorem}
 
 \begin{proof} See Matumoto~\cite{ma}.
 \end{proof}
 \vskip.1in


\section{Classes of order 4}

In this section we show that all classes of order 4 in $\calg_\qq$ in the image of $\calg^\zz$ remain of order 4 when projected to $\calg_{\qq_p}$ for some   $p  |   \Delta_V(-1)$,  $p \equiv 3 \mod 4$, where $V$ is an integer matrix representing the class in $\calg^\zz$.  We also develop simple effective criteria to detect elements of order 4 that do not require a detailed $p$--adic analysis.  The results here strengthen a result of  Morita~\cite{mo} in which it was shown that one can restrict to primes   $p$ dividing $2\Delta_V(1)\Delta_V(-1)$.  

The restriction to $p \equiv 3 \mod 4$ is automatic, given that $W(\ff_p)$ does not contain 4--torsion if $p \equiv 1 \mod 4$. The hardest technical work is in ruling out $p = 2$.

\begin{theorem}\label{main} If a class $\alpha \in \calg_\qq$ that arises from a knot $K$, or, equivalently, in the image of $\calg^\zz$, is of order 4, then $\alpha$ is of order 4 in $\calg_{\qq_p}$ for some $p \equiv 3 \mod 4$ with $p$ dividing $\Delta_V(-1)$.
\end{theorem}

The proof will use the following lemma.

\begin{lemma}\label{lemma1} Let $(M, Q, T) \in \calg_{\qq_p}^\delta$ be an isometric structure, where $p$ is odd and $\delta \in \qq_p[t^{\pm 1}]$ is monic, irreducible and symmetric.  If $\Delta_T(1)\Delta_T(-1) = p^{2e}u$ where $u$ is a unit in $\zz_p$, then $(M, Q, T)$ is not of order 4. In particular, since $\Delta_T(t) =  \delta(t)^k$,   if $\delta(1)\delta(-1) = p^{2e} u$, then $(M, Q, T)$ is not of order 4.
\end{lemma}
\begin{proof}  The form $Q$ can be diagonalized to be $[d_1, d_2,  \ldots, d_k, pd_{k+1}, \ldots ,   pd_{2d}]$, where the $d_i$ are units in $\zz_p$.  In $\qq_p^*/ (\qq_p^*)^2$, the determinant of the form is $\Delta_T(1)\Delta_T(-1) = p^{2e}u$.  Thus $k$ is even, say $k = 2l$.

  Under the isomorphism of $W(\qq_p) \cong W(\ff_p) \oplus W(\ff_p)$, $Q$ maps to the pair of forms $[d_1, d_2,  \ldots, d_{2l}] \oplus [  d_{2l+1}, \ldots ,    d_{2d}]$.  But a form of order 4 in $W(\ff_p)$ is of odd rank.  Thus, $Q$ is of order at most 2, and, applying Theorem~\ref{levine3}, $2(M, Q, T)$ is Witt trivial.

\end{proof}

 \begin{proof}[Proof (Theorem~\ref{main})] 
Let  $(M,Q, T)$ be the rational isometric structure  that arises from the Seifert matrix $V$ representing a class in $\calg^\zz$.

Fix for now a prime number $p$ that does not divide $\Delta_V(-1)$. We will show that $(M,Q,T)$ cannot represent an element of order 4 in $\calg_{\qq_p}$.

Recall that $\Delta_V(t) = \det(V)  \Delta_T(t)$.    By Gauss's Lemma, applied in the setting  $\zz_p \subset \qq_p$, we can form the $p$-adic irreducible factorization
$\Delta_V  =   \prod_i \dot{\delta}_i \prod_j \dot{f}_j$
where the $\dot{\delta}_i \in \zz_p[t]$ are the symmetric factors and the remaining factors, the  $\dot{f}_j\in \zz_p[t,t^{-1}] $,  occur in $(t \to t^{-1})$ conjugate pairs.    The dots over the polynomials indicate that these are associates (differ by multiplication by a nonzero element of $\qq_p$) of the irreducible monic factors of $\Delta_T$.

If $(M,Q, T)$ is of order 4 in $\calg_{\qq_p}$, then the image of $(M,Q, T)$ in $\calg_{\qq_p}^{\delta_i}$ will be of order 4 for one of the $\delta_i$, which we now denote $\delta$.  Call this image   $(M_\delta, Q_\delta, T_\delta)$.   \vskip.1in

\noindent{\bf Case I, $p$ odd:}  (Morita's theorem)  Since $p$ does not divide  $\Delta_V(1)\Delta_V(-1)$, this product  is a unit in $\zz_p$, and the same is true for $\dot{\delta}(1)\dot{\delta}(-1)$.  
It follows that $\delta(1)\delta(-1)$ is of the
   form $a^2 u$ where $u$ is a unit, and thus Lemma~\ref{lemma1} applies to show that $(M_\delta, Q_\delta, T_\delta)$ is not of order 4.  \vskip.1in

\noindent{\bf Case II, $p=2$:}     Recall that the {\it discriminant} of a form $Q$ over
 $\ff$ of even rank $2e$ is defined
 to be  $\disc (Q) = (-1)^{e}\det(Q)$.  (See Appendix~\ref{wittgroupsection}.)  We capitalize and use the symbol $\Disc$ to designate the discriminant of a polynomial.)  This determines
 a homomorphism $I \to \ff^*/(\ff^*)^2$,
  where $I$ is the subgroup of $W(\ff)$ generated by forms of even rank.

In the present situation,  as mentioned above,  we have from~\cite{sch} that $\qq_2^*/(\qq_2^*)^2$ is of order 8, with representatives  $\{ \pm 1 , \pm 2, \pm 5, \pm10\}$. One then checks immediately that   the values  of the discriminants of the classes of order 4  are $\{ -1 , -2, -5, -10\} \subset \ \qq_2^*/(\qq_2^*)^2$.   It remains to show that $\disc(Q_\delta)$ is not in $\{ -1 , -2, -5, -10\}$ modulo squares.

The characteristic polynomial of $T_\delta$ is, up to a constant, $g(t)$ where $g(t)$ is a symmetric polynomial with coefficients in $\zz_2$.  Since $\Delta_V(1) = 1$ and $g(t)$ is a factor of $\Delta_V(t)$ in the $\zz_p$--adic factorization,   $g(1)$ is a unit in $\zz_2$,  and so after multiplying by another constant we can assume that $g(1) = 1$.  Given that $g$ is symmetric and of even degree $2e$, we write
$$g(t) = a_0 + a_1 t + \cdots + a_{e-1}t^{e-1} + a_e t^e +  a_{e-1}t^{e+1}+ \cdots + a_0t^{2e}.$$ Since $g(1) = 1$, we have $a_e = 1 - 2a_{even} - 2a_{odd}$, where $a_{even}$ and $a_{odd}$ are the sums of the coefficients with even or odd index, respectively.

Since the determinant of $Q_\delta$ is given by $g(1)g(-1)$ modulo squares, the discrimant of $Q_\delta$ is given by $(-1)^e g(1)g(-1) = (-1)^e g(-1)$, which expanded equals 
$$(-1)^e ( 2a_{even} - 2a_{odd} + (-1)^e ( 1 - 2a_{ even} - 2 a_{odd}) =1 + 4a^*,$$  for some $a^*$.  Thus, $\disc(Q) \equiv 1 \mod 4$.  None of the  elements in $\{ -1 , -2, -5, -10\}$ are equivalent to 
$1 \mod 4$   and 
 thus $Q_\delta$ is not of order 4 in $W(\qq_2)$.
\end{proof}

\begin{corollary}\label{cor1} If $\Delta_V(-1)$ has no prime factor $p$ with $p \equiv 3 \mod 4$, then $K$ is not of order 4 in $\calg$.
\end{corollary}

\begin{corollary}\label{cor2} If $V$ is of order 4 in $\calg$ then for some $p  \equiv 3 \mod 4$ and some symmetric irreducible factor $g(t) \in \zz[t,t^{-1}]$ of $\Delta_V(t)$, $p$ divides $g(-1)$ and $g$ has odd exponent in $\Delta_V$.
\end{corollary}

\begin{proof} If $V$ has order 4,   this will be detected  in $\calg_{\qq_p}$ for some  $p \equiv 3 \mod 4$ that divides $\Delta_V(-1)$.  In turn, it will be detected in $\calg_{\qq_p}^\delta$ for some $\delta$ that divides $\Delta_V(t)$.  Suppose that $g(t)$ is the irreducible factor of $\Delta_T(t)$ that is divisible by $\delta(t)$.  If $g(t)$ has even exponent in $\Delta_T(t)$, then $\Delta_{T}^\delta$ will be an even power of $\delta$.  Thus, according to Lemma~\ref{lemma1}, $2(M^\delta, Q^\delta, T^\delta) = 0 \in \calg_{\qq_p}^\delta$.

\end{proof}



  
 


\section{Classes of Order 2}

In this section we consider forms $(M,Q,T) \in \calg_\qq$ that are known to be of finite   order and not of order 4.  There are two cases, one of which is trivial.

\subsection{The trivial case:  odd exponent}$\ $
Suppose the $\Delta_T(t)$ has a symmetric irreducible  factor with odd exponent.  Then 
 $(M,Q,T)$ is nontrivial, and so of order exactly 2.
 
 \subsection{The even exponent case}$ \ $  We are reduced to the case that $(M,Q,T)$ is of order 1 or 2, and all irreducible symmetric factors of $\Delta_T(t)$ have even exponent.
  In the case of identifying order 4 classes, we saw that primes that divide the determinant of the class were key.  Here we must also consider the discriminant of the polynomial, $\Disc(\Delta_T)$, and $\det(V)$.  The definition of these is presented in Appendix~\ref{discsection}.

  \begin{theorem} Let $V$ be a nonsingular Seifert matrix representing a class in $\calg^\zz$ of rank $2g$ and let $(\qq^{2g}, V+V^t, V^{-1}V^t) = (M,Q,T) \in \calg_\qq$.    Suppose that all irreducible symmetric factors of $\Delta_V(t)$ have even exponent.  Then for any prime $p$ that does not divide $2 \det( V) \Disc (\bar{\Delta}_V(t))$, $(M,Q,T) = 0 \in \calg_{\qq_p}$, where $\bar{\Delta}_V(t)$ denotes the product of all the distinct irreducible factors of $\Delta_T$. 
  
  \end{theorem}
 
 \begin{proof}We begin by defining $ \calg_{\zz_p}$.  This is the Witt group consisting of triples $(M,Q,T)$ where $M$ is a free $\zz_p$--module, $Q$ is a symmetric bilinear form on $M$ with determinant a unit in $\zz_p$, and $T$ is an isometry of $(M,Q)$.
 
We show in Lemma~\ref{vpluslemma} below that since $p$ does not divide $\Disc(\bar{\Delta}_V)$,  $det(V+V^t)$ is a unit in $\zz_p$.  Furthermore, the entries of  $V^{-1}V^t$ are rational with denominators prime to $p$ and so all entries are elements of $\zz_p$.  It follows that this matrix defines an isometry of $\zz_p^{2g}$.  Thus, the class $(\qq^{2g}, V+V^t, V^{-1}V^t) \in \calg_\qq$ is in the image of the class $ (\zz_p^{2g}, V+V^t, V^{-1}V^t) \in \calg_{\zz_p}$

 The characteristic polynomial $\Delta_T$ of $V^{-1}V^t$ is a monic polynomial in $\zz_p[t^{\pm1}]$, and  has irreducible factorization over $\zz_p$ (and $\qq_p$) as $\prod_i \delta_i^{\epsilon_i}\prod_j   g_j$, where the $\delta_i$ are the  irreducible symmetric factors and the $\epsilon_i$ are all assumed to be even.
 
 Since $p$ is prime to the discriminant,$\Disc(\Delta_T)$, as we describe in Appendix~\ref{discsection} the splitting $(M,Q,T) = \oplus (M^{\delta_i}, Q^{\delta_i}, T^{\delta_i})$ can be viewed as a splitting in $\calg_{\zz_p}$, rather than in $\calg_{\qq_p}$. 
 
  It remains to show that each one of these summands, say $ (M^{\delta}, Q^{\delta}, T^{\delta})$, is Witt trivial  as a class in $\calg_{\qq_p}^\delta$.  Since the exponent is even, by Theorem~\ref{levine3} it is sufficient to show that $(M^\delta, Q^\delta) = 0 \in W(\qq_p)$.  Notice that since the exponent on $\delta$ is even, $M$ has rank $4m$ for some $m$ and the determinant (and thus the discriminant) of $Q$ is a square, and so is trivial in $\qq_p^* / (\qq_p^*)^2$.
  
  Diagonalize $Q^\delta$ to be $[ u_1, \ldots , u_k, pv_1, \ldots , pv_j]$.  There is   the homomorphism $\phi^o \co W(\qq_p) \to W(\ff_p)$ that takes our diagonalized class to $[v_1, \ldots , v_j]$.  However, as described in Appendix~\ref{wittgroupsection}, $\phi^o$ vanishes on $W(\zz_p)$.  Thus, $[v_1, \ldots , v_j]$ is of even rank and has discriminant 1.  It follows that applying $\phi^e$ to our form (resulting in $ [ u_1, \ldots , u_k]$) is also a form of even rank and discriminant 1, and so is trivial in $W(\ff_p)$. 
  
  Since  $p$ is odd, $\phi^o \oplus \phi^e$ defines an isomorphism $W(\qq_p) \to W(\ff_p) \oplus W(\ff_p)$.  Thus we see that $(M^\delta, Q^\delta)$ is Witt trivial.   
 
 \end{proof}
 
 \begin{lemma}\label{vpluslemma} Let $V$ be a Seifert matrix and suppose the prime $p$ does not divide $\det(V)\Disc(\bar{\Delta}_V(t))$. Then $\det(V + V^t)$ is a unit in $\zz_p$.
 \end{lemma}
 
 \begin{proof}  We begin by noting that $\det(V + V^t) = \Delta_V(1)\Delta_V(-1)$.  This will be a unit if and only if $\bar{\Delta}_V(1)\bar{\Delta}_V(-1)$ is a unit; removing multiple factors does not change whether an element in $\zz$ is divisible by $p$.  The leading coefficient of $\bar{\Delta}_V(t)$ is a divisor of $\det(V)$, and so is prime to $p$ and is a unit in $\zz_p$.  Dividing by that leading coefficient yields a monic polynomial $\bar{\Delta}_T(t) \in \zz_p$.  We need to show that $\bar{\Delta}_T(1)\bar{\Delta}_T(-1)$ is a unit in $\zz_p$.
 
 The discriminant of $\bar{\Delta}_T(t)$ is given by $ \Disc(\bar{\Delta}_T(t)) = \prod_{i,j} (\alpha_i - \alpha_j)^2$ where the product is over all distinct pairs of roots of $\bar{\Delta}_T(t)$ in the algebraic closure of $\qq_p$.  Since $ \bar{\Delta}_T(t)$ is symmetric and does not have $\pm 1$ as  a root, if $\alpha$ is a root, then so is  $  \frac{1}{\alpha} \ne \alpha$.  Collecting roots that occur in inverse pairs, we find  
  $$\Disc(\bar{\Delta}_T(t)) =  \prod  (\alpha_i - \frac{1}{\alpha_i})^2\prod(\alpha_i - 
 \alpha_j)^2, $$ where the second product is taken over pairs with $\alpha_j \ne   \frac{1}{\alpha_i}$.

 This product can be rewritten as 
 $$ \frac{1}{\prod  {\alpha_i^2}} \prod_i (\alpha_i^2 - 1)^2\prod(\alpha_i - 
 \alpha_j)^2   =  \frac{1}{\prod_i {\alpha_i^2}} \prod  (\alpha_i - 1)^2(\alpha_i +1)^2\prod(\alpha_i - 
 \alpha_j)^2. $$
 Since $\bar{\Delta}_T(t)$ is monic and symmetric, $\prod \alpha_i =1$ and the entire product can be simplified to give 
 $$\Disc(\bar{\Delta}_T(t)) =  \bar{\Delta}_T(1)^2\bar{\Delta}_T(-1)^2\prod(\alpha_i - 
 \alpha_j)^2.$$  The two elements $ \bar{\Delta}_T(1)$ and $\bar{\Delta}_T(-1)$ are clearly in $\zz_p$, and we are assuming that $ \Disc(\bar{\Delta}_T(t))$ is a unit in $\zz_p$.  Thus, if we show that $ \prod(\alpha_i - 
 \alpha_j)^2$ is in $\zz_p$, each of $ \bar{\Delta}_T(1)$,  $\bar{\Delta}_T(-1)$, and    $ \prod(\alpha_i - 
 \alpha_j)^2$ are seen to be units in $\zz_p$.  
 
 To see that $ \prod(\alpha_i - 
 \alpha_j)^2$ is in $\zz_p$, note that it is fixed by the Galois group of the splitting field of $\bar{\Delta}_T(t)$ over $\qq_p$, and thus is in $\qq_p$.  However, it is an algebraic integer in the algebraic closure of the fraction field of $\zz_p$.  The only algebraic integers in the fraction field of an integral domain must be in the domain itself.  Thus, it is in $\zz_p$ as desired.

 \end{proof}

\section{The algebraic order of prime knots with 12 or fewer crossings}

There are 2,977 prime knots with      12 or fewer crossings.  Here we describe how the algebraic orders of all such knots are determined. Our goal is to present enough of  the calculation to illustrate how the complete set of results,  appearing in  the {\it KnotInfo} table~\cite{lc}, were derived.  In addition, we have isolated out special cases to demonstrate methods that readily apply when specific theorems cannot be quoted directly.

\vskip.1in
\noindent{\bf Infinite order:} If a knot has infinite order, it is detected by a nontrivial signature.  For 2,132 of the knots, the signature of $V + V^t$ is nonzero.  Another 125 knots have nontrivial $\omega$--signature (the signature of $(1-\omega)V + (1-\omega^{-1})V^t)$ nontrivial for some unit complex number $\omega$.
This leaves 720 knots of finite algebraic order.\vskip.1in

\noindent{\bf Slice knots:} There are 157 knots that have been identified as topologically slice (many through the unpublished work of Stoimenow~\cite{sto}).  This leaves 563 knots to resolve.\vskip.1in

\vskip.1in
\noindent{\bf Order 4:}  By Theorem~\ref{oddexp}, if $K$ is of finite algebraic order and $D =\Delta_K(-1) =  \det(V + V^t)$ is divisible by a prime $p$, $ p \equiv 3 \mod 4$ and $p$ has odd exponent in $D$, then $K$ is of order 4. This applies to 172 of the remaining knots, leaving 391 to resolve.
\vskip.1in

\noindent{\bf Order 2:} If $K$ is of finite algebraic order and its Alexander polynomial has a symmetric factor (over $\qq$) with odd exponent, it has order 2 or 4.  If no prime $p \equiv 3 \mod 4$ divides $\Delta_K(-1)$ then $K$ is of order 2 in $\calg$.  This applies to 318 of the remaining knots, leaving 73 to resolve.

\vskip.1in
\noindent{\bf More Order 2:}  As in the previous paragraph, if $K$ has finite algebraic order and its  Alexander polynomial has a symmetric factor (over $\qq$) with odd exponent, it has order 2 or 4.  If in addition $K$ is amphicheiral (equal to its mirror image, regardless of orientation) it is of algebraic order exactly 2.  (Reversing the orientation of a knot has the effect of transposing the Seifert matrix $V$. An examination of Levine's classification reveals the this does not change the algebraic concordance class of a knot.  As a more explicit proof, see for instance~\cite{kl}.) This applies to another 5 knots, leaving 68 to resolve.
\vskip.1in

\noindent{\bf Basic examples of  algebraically slice knots:}  If the Alexander polynomial has no irreducible symmetric factors, the knot is algebraically slice.  This applies to another 9 knots, leaving 59 cases to resolve.
\vskip.1in

\noindent{\bf More Order 2:}  Of these remaining 59 knots, 25 have the property that an irreducible factor $g(t)$ of the Alexander polynomial has odd exponent, so that the knot is of order 2 or 4, but no such irreducible factor has $g(-1)$ divisible by a prime $p \equiv 3 \mod 4$.  Thus, by Corollary~\ref{cor2}, the knot has order exactly 2. 
As an example, the knot $9_{24}$ has Alexander polynomial 	$1-5t+10t^2-13t^3+10t^4-5t^5+t^6$ with determinant 45 = $3^2 5$. Thus it is conceivable that it could be of order 4, detected at the prime 3.  But the polynomial factors as $(1 - 3 t + t^2) (1 - t + t^2)^2$, and the 3 factor arises from a symmetric irreducible factor of exponent 2.  As a second example, the knot $9_{34}$ has Alexander polynomial $   $ with determinant $3^2 5$. In this case the polynomial factors as $(-2 + t) (-1 + 2 t) (1 - 3 t + t^2)$ and thus the 3 factor does not arise from a symmetric factor of the Alexander polynomial. All 25 cases are similar to one of these two examples.  

This leaves 34 knots to consider.

\vskip.1in

\noindent{\bf More Algebraically Slice:} There are  nine remaining knots for which all symmetric factors have even degree.  These  are either trivial in $\calg_\qq$ or are of order 2.  We can rule out order 2 in all the cases to see that these are algebraically slice, as follows.

For seven of these knots, there is a unique symmetric irreducible factor $\delta$, and it is of degree 2.  For instance, for $12_{a169}$ the Alexander polynomial factors as $(2 - 3t + 2t^2)^2$ and for $12_{n224}$ the Alexander polynomial factors as $(1-2t)(2-t)(1-t+t^2)^2$.  
In this case, regardless of the prime $p$, if the quadratic factors over $\qq_p$, then the form is Witt trivial, since no even degree symmetric factors would remain.  If the quadratic is irreducible, then the form $Q = V + V^t$ would be Witt trivial if and only if the form $Q_\delta$ is Witt trivial, since the form $Q$ restricted to the complement of the $\delta$ summand is automatically Witt trivial.  In each case, one can diagonalize the form $V +V^t$ and find that the diagonal elements are paired $[d_1, -d_1, d_2, -d_2, \ldots]$.

One of the two more challenging cases is that of $K = 12_{a990}$.  There is a nonsingular Seifert matrix $V$ for $K$ of size $8 \times 8$, and  $\Delta_K(t) = (t^2 -t +1)^2(t^2 -3t +1)^2$.  Thus the corresponding transformation $T$ has characteristic polynomial 
$ (t^2 +t +1)^2(t^2 +3t +1)^2$.  Letting $\delta_1(t) = t^2 +t +1$ and $\delta_2(t) =  t^2 +3t +1 $ it is easy to find a basis for $M_{\delta_1}$ and  $M_{\delta_2}$ over the rationals.  These are just the images of the transformations $\delta_2(T)^2$ and $\delta_1(T)^2$ respectively, both of which are rank 4.  On each of these, it is easy to find the respective quadratic form, simply by restricting $Q$ to each, and these can be diagonalized over the rationals, with diagonal entries integers.  For some primes $p$,  $\delta_i$ might factor over $\qq_p$, but since it is quadratic, if it factors then the form is automatically Witt trivial.  So we assume that $\delta_i$ is irreducible over $\qq_p$.  In this case, it is sufficient to show that the forms $Q_{\delta_i}$ are trivial over the rationals (the exponent of the characteristic polynomial is even).  For this, one can apply the classification of Witt forms over $\qq$, as described in Appendix~\ref{wittgroupsection}.   This calls for a consideration of all prime integers, but if $p$ does not divide any of the diagonal entries of $Q_{\delta_i}$, then one notes that the induced forms in $W(\ff_p)$ are of rank 4 with trivial discriminant (that is, a square), and thus the forms are Witt trivial.  At the finite set of primes that remain, one must check that the image forms in $W(\ff_p)$ are trivial.  In the actual calculation for these examples, only four primes appear (though this might depend on the choice of spanning set of $M_{\delta_i}$) and so the calculation is quickly done.

The second case that requires further calculation is that of $12_{n681}$ which has a nonsingular Seifert matrix of size $8 \times 8$ and Alexander polynomial $(t^4 - t^3 + t^2 -t +1)^2$, so that the corresponding transformation $T$ has characteristic polynomial  $(t^4 - t^3 + t^2  - t +1)^2$.  Letting $\delta(t) =    t^4 - t^3 + t^2 -t +1$, one finds the image of the transformation $\delta(T)$ is a rank 4 invariant subspace of the 8--dimensional rational vector space on which the form vanishes.

\vskip.1in
\noindent{\bf Generalizing this example} 

 In the case that  the Alexander polynomial 
  of $K$ factors as $\delta(t)^2$ with $\delta(t)$ irreducible, if the transformation $\delta(T)$ is nontrivial then $K$ is algebraically slice.  The proof consists of
  noting that the $\qq[t^{\pm1}]$--module
    $M$ is 
  of the form  $\qq[t^{\pm1}] / <\delta(t)^2>$. 
   The image of $\delta(T)$ will be half dimensional 
  and the form $Q$ vanishes on the image.  In the previous example, that of $K_{a990}$, this approach does not work:  as a $\qq[t^{\pm1}]$--module 
 $M$ is isomorphic to  $\qq[t^{\pm1}] / <(t^2 + 3t + 1)^2> \oplus  \qq[t^{\pm1}] / <t^2 +  t + 1>\oplus  \qq[t^{\pm1}] / <t^2 +t   + 1> $ and thus there are an infinite number of invariant submodules to consider.
\vskip.1in

\noindent{\bf Order 2:}  The  remaining  25 cases are the most technical.  In these cases there is a unique prime $p \equiv 3 \mod 4$ dividing $\Delta_K(-1)$ and in each case $p$ has exponent $2$ in $\Delta_K(-1)$.  Furthermore, there is an irreducible factor $f$ of $\Delta_K$ (over $\zz$) such that $f$ has exponent $1$ in $\Delta_K$ and  $p^2 $ divides $f(-1)$.  If $f$ remains irreducible in $\qq_p$ then arguments as given above would imply that $K$ is of order  exactly 2.

However, it is conceivable that $K$ will be of order 4, but for this to occur,   $f$  would factor as $f = f_1 f_2 \cdots f_n$ over $\qq_p$ (and thus over $\zz_p$), and for at least one of the symmetric $f_i$ (and so for at least two of the symmetric $f_i$), $ f_i(1)f_i(-1)  \ne 1 \in \qq_p^* / (\qq_p^*)^2$, or, put otherwise,  $D(f_i(1)f_i(-1)) \ne 1 \in \zz/2\zz$ (where $D \co \qq_p^* / (\qq_p^*)^2 \to \zz/2\zz$ is defined in Appendix~\ref{padicsection}).  Call two of these factors $f_a$ and $f_b$.

 There is the canonical homomorphism $\zz_p \to \zz / p\zz$. Denote this map by $x \to \bar{x}$.  Similarly, there is the induced map $\zz_p[t] \to (\zz / p \zz)[t]$, which we denote by $f(x) \to \bar{f}(x)$.

 The symmetric factors $f_a$ and $f_b$ are both of even degree, and thus $\bar{f_a}$ and $\bar{f}_b$ are even degree and symmetric, after factoring out a power of $t$ so that they have nonzero constant term.

 Since $D(f_a(1)f_a(-1)) \ne 1$ (that is, $f(1)f(-1) = p^k u$ where $k$ is odd and $u$ is a unit), it must be the case that $\bar{f_a}(1)\bar{f_a}(-1) = 0 \in \zz / p\zz$.  Thus, $\bar{f_a}$ must be divisible by $t \pm 1$.  Similarly for $\bar{f_b}$. We can show that each of the 25 knots in this category don't satisfy this criteria.  A few examples follow.
 
 The knot $11_{a300}$ has $\Delta_K(1)\Delta_K(-1) = 3^2 17$.  Thus the only prime of interest is $3$.  When reduced modulo 3, we have the irreducible factorization $\bar{\Delta}_K(t) = (1 + t)^2 (1 + t^2) (1 + t + t^2 + t^3 + t^4)$.  These factors cannot be combined to find two symmetric factors, each of which is of even degree and divisible by $ t \pm 1$.
 
 A similar example is $12_{a1170}$, again with $\Delta_K(1)\Delta_K(-1) = 3^2 17$.  Its Alexander polynomial reduced modulo 3 satisfies  $\bar{\Delta}_K(t) = 2 (1 + t)^2 (2 + t + t^2 + t^3) (2 + 2 t + 2 t^2 + t^3)$.  In this case, by distributing the $1+t$ factors between the two other factors we split the polynomial into  even degree polynomials, each of which evaluates trivially at $t = -1$.  However, neither of these is symmetric.
 
 This approach works for 24 of the 25 knots of this variety.  The one exception is $12_{n525}$.  It has $\Delta_K(t) = 1-8t+28t^2-43t^3+28t^4-8t^5+t^6$.    Again, $\Delta_K(1)\Delta_K(-1) = 3^2 17$.

 Working modulo 3, this polynomial factors as $(1+t)^4 (1+t^2)$. Suppose that $\Delta_K(t)$ factors nontrivially with two or more symmetric factors, each of even degree, over the $3$--adics.  One possibility would be that there are  degree 2 and degree 4 irreducible factors.  In this case, one possibility for the corresponding factorization modulo 3 would be $[(1+t ^2)][(1+t)^4]$ and the other would be $[(1+t)^2][(1+t)^2(1+t^2)]$.  
The other possibility is that $\Delta_K(t)$ factors over the $3$--adics as the product of three symmetric quadratics.  Then, modulo 3, the corresponding factorization would be $[(1+t)^2] [(1+t)^2] [(1+t^2)]$.
 
 Notice that among these three cases, there are only two cases in which there are two irreducible factors over the $3$--adics both of which  satisfy $\delta(1)\delta(-1) = 0 \mod 3$.  In each of these two cases there is a quadratic factor which satisfies $\delta(1)\delta(-1) = 0 \mod 3$.  We want to show that this does not occur.

 One way to do this is to find the $p$--adic factorization.  Another way is to check for factorizations modulo $3^k$ for various $k$.  The second method can be done quickly by computer, and we find that modulo 27 the only factorization into a quadratic and quartic has quadratic term $1 + t^2 \mod 3$, and this does not satisfy $\delta(1)\delta(-1) = 0 \mod 3$.
 
 We  now want to describe a method for factoring over the $p$--adics.  
 If $\Delta_K$ factored as desired, then we would have  $\Delta_K = f g$, where $f(t) = 1 + at +t^2$ and $g(t) = 1 + bt + ct^2 + bt^3 + t^4$, where $a, b , c \in \zz_3$.  (This uses Gauss's Lemma and the fact the $\Delta_K(t)$ is monic.) 
 
 For this to hold, we see immediately that $b = -8 -a$.  Making this substitution into $g$ and multiplying gives $ fg - \Delta_K = (-27-8 a-a^2+c)t^4+(27-2 a+c a) t^3+(-27-8 a-a^2+c)t^2$, so $c = 27 + 8a + a^2$.
 
 Again substituting and expanding gives  $ fg - \Delta_K =  (27+ 25a+8a^2+a^3)t^3$.  The number $a = 0$ is a   solution modulo 3 to the equation $h(a) = 27 + 25a+8a^2+a^3 =0$. According to the general theory of $p$-adic polynomials, since $0$ is not a solution to $h'(a) = 0$, it lifts to a $p$-adic solution.  In this case, it is relatively easy to find that lifting: knowing the value mod 3 permits one to find the solution modulo 9; this solution then is easily lifted to a solution modulo 27, and so on.  For instance,    modulo $3^8$ a solution is $a = 2565 = 2(3^3) + 1(3^4)+1(3^5) +1(3^7)$.  The factorization of $\Delta_K$ modulo $3^8 = 656144 $ is $$\Delta_K(t) = (1+2565 t+t^2) (1+3988 t + 5967 t^2+3988 t^3 +t^4) \mod 3^8 .$$

\appendix

\section{Background: $p$-adic numbers}\label{padicsection}

A good concise introduction to $p$--adic arithmetic is contained in the text by Serre~\cite{se}.
Let $p$ be a prime integer.  A $p$--adic integer can be defined to be a formal sum   $\sum_{i=0} ^\infty a_i p^i$, where the   $a_i$ satisfy $0 \le a_i < p$.  The set of $p$--adic integers is denoted $\zz_p$. Addition and multiplication are defined formally, and with these operations $\zz_p$ forms a commutative ring with unity.  There is a natural inclusion of $\zz$ into $\zz_p$.  The set of units in $\zz_p$ are those numbers for which $a_0 \ne 0$.  Up to multiplication by a unit, $p$ is the unique prime in $\zz_p$.

The $p$--adic rationals, $\qq_p$, are defined to be formal sums $\sum_{i=k} ^\infty a_i p^i$ where   $a_i$ satisfies $0 \le a_i < p$, but we no longer assume that $k \ge 0$.  Since any element in $\qq_p$ can be written as $p^i a$ where $a$ is a unit in $\zz_p$, it is clear  that $\qq_p$ is a field, and since it is the minimal field that contains $1/p$, it is the field for fractions of $\zz_p$.

The following function $D$ will be of use.  Let $\qq_p^*$ denote the nonzero elements in $\qq_p$.

\begin{definition}  Any nonzero element $x \in \qq_p^*$ can be written as $p^\epsilon u$ where $u$ is a unit in $\zz_p$.  Set $D(x) \in \zz/2\zz$ to be the mod 2 reduction of $\epsilon$.

\end{definition}
\vskip.1in
\noindent{\bf The structure of $\qq_p^*/ (\qq_p^*)^2$}

\vskip.1in

 \begin{theorem}   For $p$ odd, the quotient of the multiplicative group of $\qq_p$  by its squares,   $\qq_p^*/ (\qq_p^*)^2$, is isomorphic to $\zz/2\zz \oplus \zz/2\zz $.  Four distinct elements are given by the set  $S = \{1, u,   p,     up\} \subset \qq_p^*/ (\qq_p^*)^2$, where $u$ is any  integer $0< u < p$ which is not a square modulo $p$. 
 \end{theorem}

 \vskip.1in
  \begin{proof}
Let $\ff_p$ denote the field with $p$ elements, $\zz /  p\zz $.  Its multiplicative subgroup, $\ff_p^*$, is a cyclic group of even order $p-1$.  It follows that  $\ff_p^*  / (\ff_p^*)^2  = \zz/2\zz$.  Thus, there is an element $u$ that is not a square.

 Clearly all elements of $\qq_p^*/ (\qq_p^*)^2$ are of order 2.  It is easily seen that   $S$ is a subgroup and that no product of any pair of distinct elements in $S$  is a square, so $S$ is a subgroup isomorphic to  $\zz/2\zz \oplus \zz/2\zz $.

 Any element in $\qq_p^*$ can be multiplied by an even power of $p$ so that it is of the form $ s (a_0 + a_1p + a_2p^2 + \cdots)$, where $s \in S$ and $a_0$ is a square module $p$.  Finally, a square root to $ (a_0 +  a_1p + a_2p^2 + \cdots)$ is easily found, using the fact the $a_0$ is a square modulo $p$ and solving recursively for the coefficients.\vskip.1in

The case of $p = 2$ is a bit more delicate, and we leave the proof to~\cite{sch}.\vskip.1in
\end{proof}

 \begin{theorem}   The quotient of the multiplicative group of $\qq_2$  by its squares,   $\qq_2^*/ (\qq_2^*)^2$, is isomorphic to $\zz/2\zz \oplus \zz/2\zz \oplus \zz/2\zz$.  Eight distinct elements are given by the set  $S = \{\pm 1, \pm 2,  \pm 5, \pm 10     \} \subset \qq_2^*/ (\qq_2^*)^2$.
 \end{theorem}

\vskip.1in

Finally we note that the function $D$, defined   earlier, descends to a function on $\qq_p^* /(\qq_p^*)^2$, which we denote $\bar{D}$.

\section{Background: Witt groups}\label{wittgroupsection}

The basic theory of Witt groups of symmetric bilinear forms can be found in~\cite{mh, sch}.   For the most part we are interested only in forms over fields.  In one case we need to consider the ring $R = \zz_p$, so we state the basic definitions in terms of a more general commutative ring $R$ with unity and finitely generated free modules over $R$ instead of vector spaces. To be more specific,  let $R$ be a ring, either the field $\qq$, $\rr$, $\qq_p$, or $\ff_p$, the finite field with $p$ elements, where $p$ is prime, or the ring $\zz$ or $\zz_p$.

Consider pairs $(M,Q)$ where $M$ is a free module over $R$ and $Q$ is a nonsingular  symmetric bilinear form on $M$.  Here nonsingular means that the determinant of $Q$ is a unit.    The form $(M,Q)$ is called Witt trivial if $M$ is of dimension $2g$ for some $g$ and there is a summand of $M$ of dimension $g$ on which $Q$ is identically 0.  Such a summand is called a metabolizer for $(M,Q)$.   Forms $(M_1, Q_1)$ and $(M_2, Q_2)$ are called Witt equivalent if $(M_1 \oplus M_2, Q_1 \oplus -Q_2)$ is Witt trivial.

 The set of Witt equivalence classes  of pairs $(M,Q)$ constitute an  abelian group under the operation induced by direct sum, and this is the Witt group  $ W(R)$.

If $R$ is not of characteristic 2 (that is, all rings under consideration except $\ff_2$),  any form $Q$ has a diagonal matrix representation with respect to some basis of $M$. In the case of $\ff_2$, not every form is diagonalizable,  but it is true  that all Witt classes are represented by diagonal forms.  We abbreviate such a  diagonalization with the vector of its diagonal entries:  $[d_1, \ldots , d_k]$.  Notice that if $Q$ is diagonalized with respect to some basis and if a basis element is replaced with a multiple, then the corresponding diagonal entry is multiplied by the square of that constant.

There are two basic functions defined on $W(R)$. 

\begin{definition} The  rank of a class $w \in W(R)$ is defined to be the rank of a representative of $w$, reduced modulo 2, denoted rk$(w) \in \zz/2\zz$.  
\end{definition}

\begin{definition}  The kernel of the rank function is called the fundamental ideal, $I(\ff)$.  \end{definition}

\begin{definition} The discriminant of a class in   $w \in W(R)$ represented by a form $Q$ of rank $r$  is $\disc( w) =  (-1)^{r(r-1)/2}\det(Q)$.
\end{definition}

The discriminant is {\it not} a homomorphism on $W(R)$.  However, it is when restricted to $I(\ff)$, that is, to even rank forms.

\begin{theorem} If $Q$ is a form of rank $2g$, then $\disc( Q) = (-1)^g \det(Q) $ and  $\disc \co I(R) \to \zz / 2\zz$ is a homomorphism.  

\end{theorem}
 \vskip.1in
 \noindent{\bf Examples} 
 
 We state the following results leaving most details to the references.

 \subsection{$W(\rr) \cong \zz$}$ \ \ $

   Any Witt class $\alpha$ has a diagonal representative   $[1, \ldots, 1, -1, \ldots, -1]$.  The sum of the entries is   the signature, $\sigma(\alpha) \in \zz$.  This induces an isomorphism $\sigma: W(\rr) \to \zz$.

 \vskip.1in

 \subsection{${W(\ff_p)} = \zz/2\zz, \zz/2\zz \oplus \zz/2\zz\text{\  or\ }\zz/4\zz$, depending on $p = 2,  1,\text{\ or\ }3$    mod $  4$, respectively.}$ \ \ $

 If $p = 2$, then simultaneous row and column operations can reduce the form to a direct sum of the forms represented by the matrices  $
 \begin{pmatrix}{ }
 1
 \end{pmatrix}
$ and
$
 \begin{pmatrix}{}
  0&   1    \\
 1& 1
\end{pmatrix}.
$  The first has order 2 in $W(\ff_2)$ and the second is Witt trivial.

If $p$ is odd, let $a$ be a nonsquare mod $p$.  Then any form is equivalent to a diagonal form $[1, \ldots, 1, a, \ldots a]$.  If $p = 1 \mod 4$ then $-1$ is a square, and any diagonal  form $[b,b]$ is equivalent to $[b,-b]$, which is Witt trivial, and hence all nontrivial forms are of order 2.  Thus, every class is either trivial or one of $[1], [a]$, or $[1,a]$.  Finally, it is easily checked that none of these are Witt equivalent.

 If $p=3$ mod 4, then $-1$ is not a square, so $a$ can be taken to be $-1$.  The form $[1,-1]$ is   Witt trivial, so every class is equivalent to a positive multiple of the  diagonal form $[1]$ or the diagonal form $[-1]$.  Any form $[b,b]$ is nontrivial, again since $-1$ is not a square mod $p$, but the form $[b,b,b,b]$ is trivial, with
 metabolizer $< (1,0,a,b), (0,1,b,-a)> \in \ff_p^4$, where $(a,b)$ satisfy $1 + a^2 +b^2 = 0$.   (As described in~\cite{mh}, the existence of such an $(a,b)$ is implied by  {\it Shoebox principle} .  The set of values of $x^2$ contains $(p+1)/2$ values in $\ff_p$.  Similarly, the set of values of $1 - y^2$ contains  $(p+1)/2$ values.  Since the there are only $p$ elements in $\ff_p$, for some $x$ and $y$, $x^2 = 1 - y^2$ has a solution.)

 \begin{theorem}  A class in $W(\ff_p)$ with $p$ odd is uniquely determined by its mod 2 rank and discriminant.  For $p =2$ it is determined by its mod 2 rank.
 \end{theorem}
 
 \vskip.1in
 \subsection{For $p$ odd, $W(\qq_p) \cong W(\ff_p) \times W(\ff_p)$}$ \ \ $

  A form $Q$ over $\qq_p$ can be diagonalized as $[u_1, \ldots , u_k, pv_1, \ldots, pv_j]$ where the $u_i$ and $v_i$ are units.  If for a unit $u = a_0 + a_1p + \cdots$ we let $\bar{u}$ denote the nonzero element $a_0 \in \ff_p^*$, then we extract two forms in $W(\ff_p)$:  $[\bar{u_1}, \ldots , \bar{u_k}]$ and $[\bar{v_1}, \ldots , \bar{v_j}]$.  This map provides the desired isomorphism.

  Denote the isomorphism just defined by $\psi_p^e \oplus \psi_p^o$.
  
  \begin{theorem} There is an exact sequence $$0 \to W(\zz_p) \to W(\qq_p) \overset{\psi^o}{\to} W(\ff_p) \to 0.$$
   \end{theorem}
   This is essentially Corollary 3.3 in Chapter 4 of~\cite{mh}.  The result there applies in the more general setting in which  the last map need not be surjective. In our case surjectivity is clear.

\vskip.1in
 \subsection{For $p =2$, $W(\qq_2) \cong \zz/8\zz \oplus  \zz/2\zz \oplus \zz/2\zz$}. Here the generators are $[1], [-1,5], \text{and} [-1,2]$.  (See~\cite[Chapter 5, Theorem 6.6]{sch}.)
 \subsection{$W(\qq) \to \oplus_p W(\ff_p)$}$ \ \ $

 Any form in $Q \in W(\qq)$ can be diagonalized so that the diagonal entries are square free integers.  Fix a prime $p$ and write the diagonalized form as $$[d_1, \ldots , d_k, pd_{k+1}, \ldots , pd_{n}],$$ where the $d_i$ are all relatively prime to $p$.

The map $\psi_p^e : W(\qq) \to W(\ff_p)$ sends $Q$ to $[d_{k+1}, \ldots , d_n]$.  Combining these over all primes $p$ gives a homomorphism $W(\qq) \to \oplus_p W(\ff_p)$.   See~\cite{mh} for a proof that the kernel of this homomorphism is $W(\zz)$.

\begin{theorem}\label{oddexp} If $(M,Q) \in W(\qq)$,  $p$ is a prime  with $p \equiv 3 \mod 4$,    and det$(Q) = p^r \frac{a}{b}$ with $a$ and $b$ relatively prime to $p$ and $r$ odd, then  $\psi_p(Q)$ has order 4 in $W(\ff_p)$.
 \end{theorem}

\begin{proof}
 $\psi_p(Q)$ has odd rank.  According to the analysis of $W(\ff_p)$ for $p \equiv 3 \mod 4$ given above, if a form has odd rank, it is of order 4.
 \end{proof}


\section{Background: Discriminants, Resultants and Decompositions of Modules}\label{discsection}

Details regarding discriminants and resultants can be found in basic algebra texts, for instance~\cite{dummit}.
For a monic polynomial $p =  t^n + \cdots + a_0  \in \ff[t]$ the discriminant  is defined to be 
$$\Disc(p) =  \prod_{1\le i < j \le n}  (\alpha_i - \alpha_j)^2,$$ where the $\alpha_i$ form a complete set of roots of $p$ in the algebraic closure of $\ff$.  (Thus, if $p$ has multiple roots, $\Disc(p ) = 0$.)  

Given a second monic polynomial $q(t) =   t^m + \cdots + b_0 \in \ff[t]$, the resultant of the polynomials is defined to be
$$\Res(p,q) =  \prod_{1\le i \le n, 1 \le j \le m} (\alpha_i - \beta_j), $$ where the $\beta_j$ are a complete set of roots of $q$.

The discriminant and resultant have explicit descriptions as (integer) polynomials in the coefficients or $p$ and $q$.  Thus, if $p, q \in \dd[t]$ where $\dd$ is an principal ideal domain with field of fractions $\ff$, then $\Disc(p) \in \dd$ and $\Res(p,q) \in \dd$.  For the This can also be seen by noting that each is fixed by the appropriate Galois group, so is in $\ff$, and is an algebraic integer, so is in $\dd$.  
\begin{theorem} If $p, q \in \dd[t]$ are distinct irreducible polynomials then there are polynomials $a, b \in \dd[t]$ such that $ap + bq = \Res (p,q)$.
\end{theorem}

As a corollary, there is the following result.

\begin{corollary}Suppose that $T$ is an automorphism of $\dd^n$ with characteristic polynomial $\Delta_T $ having irreducible factorization $\Delta_T =  \prod g_i^{\epsilon_i}$.  If $\ Res(g_i, g_j) $ is a unit for all $i \ne j$, then $D = \oplus D^{g_i}$, where $ D^{g_i}$ is invariant under $T$ and  $T$ restricted to $D^{f_i}$ has characteristic polynomial $f_i^{\epsilon_i}$.
\end{corollary}

To apply this result, it is easier to work with a single discriminant rather than all the resultants.  In fact, we will be working with polynomials $\Delta \in \qq[t]$ and considering perhaps unknown factorizations in $\qq_p[t]$.

\begin{lemma}If $p, q$ are distinct irreducible monic polynomials in $\dd[t]$ then $\ Res(p,q)$ divides $\ Disc(pq)$.  In particular, if $\Disc(pq)$ is a unit in $\dd$, then so is $\Res(p,q)$.
\end{lemma}



\end{document}